\documentclass{article}
\usepackage{amsmath,amsthm,amssymb,authblk}

\newtheorem{theorem}{Theorem}[section]
\newtheorem{corollary}[theorem]{Corollary}

\newtheorem{lemma}[theorem]{Lemma}
\newtheorem{definition}[theorem]{Definition}
\newtheorem{example}[theorem]{Example}

\newtheorem{remark}[theorem]{Remark}
\numberwithin{equation}{section}
\def\pf{{\bf Proof.}~ }

\begin{document}
\title{Double Circulant Matrices}
\author{Yun Fan,\quad Hualu Liu\\
\small School of Mathematics and Statistics\\
\small Central China Normal University\\
\small Wuhan 430079, China}
\date{}
\maketitle

\insert\footins{\footnotesize{\it Email address}:
yfan@mail.ccnu.edu.cn (Yun Fan). hwlulu@aliyun.com (Hualu Liu).}

\begin{abstract}

Double circulant matrices are introduced and studied.
A formula to compute the rank $r$ of a double circulant matrix 
is exhibited; and it is shown that
any consecutive $r$ rows of the double circulant matrix
are linearly independent.
As a generalization, multiple circulant matrices are also introduced.
Two questions on square double circulant matrices are suggested.

\medskip
{\bf MSC classes:}~  15B05, 15B33, 94B15.

\end{abstract}

\section{Introduction}

Circulant matrices are an important class of matrices,
and extensively applied to numerical analysis,
cryptography, coding theory, etc.; cf. \cite{DP}, \cite{W}.

In this paper we introduce {\em double circulant matrices}
and study their properties.
This work is originally motivated by a research \cite{FH}
on quasi-cyclic codes of fractional index,
where we have a double circulant matrix, 
and we need to know how large is the rank of the matrix,
and which of the rows of the matrix form an independent set
of cardinality equal to the rank.

A circulant matrix is a square matrix which is fully specified
by one vector $(g_0,g_1,\cdots,g_{n-1})$ of length $n$,
or correspondingly, by one polynomial
$g(X)=\sum_{j=0}^{n-1}g_jX^j$ with degree $\deg g(X)<n$:
the first row is the specified vector, and
each next row is obtained by circularly shifting the previous row.
To describe double circulant matrices,
in Section \ref{g-c matrix} we study so-called
{\em generalized circulant matrices},
which are constructed with polynomials 
similarly to the usual circulant matrices,
except for that they are not necessarily square,
hence have one more parameter $m$ to specify the number of their rows.
We exhibit a formula to determine the rank $r$ of
a generalized circulant matrix, and prove that any consecutive
$r$ rows of the matrix are linearly independent
(Theorem \ref{thm for circulant} and its corollary below).

In Section \ref{d-c matrix}, we define a double circulant matrix by
the concatenation side-by-side of two generalized circulant matrices
(Definition \ref{def d-circulant} below).
Thus a double circulant matrix is parameterized by two polynomials
$g(X)=\sum_{j=0}^{n-1}g_jX^j$ and $g'(X)=\sum_{j'=0}^{n'-1}g'_{j'}X^{j'}$
with $\deg g(X)<n$ and $\deg g'(X)<n'$, and a positive integer $m$
which specifies the number of rows of the matrix.
Based on the results on generalized circulant matrices,
we obtain a formula to compute the rank $r$ of a double circulant matrix
by the parameters $g(X)$, $g'(X)$ and $m$, and also show that any consecutive
$r$ rows of the double circulant matrix are linearly independent
(Theorem~\ref{thm d-circulant} and its corollaries below).

In Section \ref{coding theory}, we apply the results on
double circulant matrices to two questions in coding theory.
The materials illustrate why
we are concerned with double circulant matrices.

In Section \ref{multiple circulant},
we extend double circulant matrices to multiple circulant matrices,
and suggest two open questions on square double circulant matrices.

\section{Generalized circulant  matrices}\label{g-c matrix}

In this section we always assume:

$\bullet$~ $m$, $n$ are positive integers;

\hangindent10mm
$\bullet$~ $F$ is a field whose characteristic ${\rm char} F$
is zero or coprime to $n$;

\hangindent10mm
$\bullet$~ $\omega$ is a primitive $n$-th root of unity
(which exists in a suitable extension of $F$ by the above assumption),
hence, $1=\omega^0$, $\omega$, $\cdots$, $\omega^{n-1}$ are
different from each other, they are just all $n$-th roots of unity.

By $F^n$ we denote the vector space consisting of sequences
$(a_0,a_1,\cdots,a_{n-1})$ with all $a_i\in F$.
For any vector $(g_0,g_1,\cdots,g_{n-1})\in F^n$, or equivalently,
for any polynomial $g(X)=\sum_{j=0}^{n-1}g_jX^j$ over $F$
with degree $\deg g(X)<n$, the matrix over $F$:
$$
\begin{pmatrix}g_0 & g_1 & \cdots & g_{n-2} & g_{n-1} \\
 g_{n-1} & g_0 &\cdots & g_{n-3}& g_{n-2} \\
 \cdots & \cdots & \cdots & \cdots & \cdots\\
 g_{2} & g_{3} & \cdots & g_0& g_{1}\\
 g_{1}  & g_{2} & \cdots & g_{n-1} & g_0
\end{pmatrix}
$$
is known as the {\em circulant matrix} associated with the polynomial $g(X)$.
Let $\cal S$ be the circularly shift operator of $F^n$, i.e.,
$${\cal S}(a_0,a_1,\cdots,a_{n-1})=(a_{n-1},a_0,\cdots,a_{n-2}),\qquad
 \forall~ (a_0,a_1,\cdots,a_{n-1})\in F^n.
$$
Then the $i$-th row (we are counting the rows from $0$) of the
above circulant matrix is just the vector
${\cal S}^i(g_0,g_1,\cdots,g_{n-1})$.

More generally, we define $M(g)_{m\times n}$ to be
the $m\times n$ matrix over $F$ whose $i$-th row
for $i=0,1,\cdots,m-1$ is the vector
${\cal S}^i(g_0,g_1,\cdots,g_{n-1})$; i.e.,
\begin{equation}\label{GCM}
M(g)_{m\times n}=\begin{pmatrix}g_0 & g_1 & \cdots & g_{n-2} & g_{n-1} \\
 g_{n-1} & g_0 &\cdots & g_{n-3}& g_{n-2} \\
 \cdots & \cdots & \cdots & \cdots & \cdots\\
 g_{n-m+2} & g_{n-m+3} & \cdots & g_{n-m}& g_{n-m+1}\\
 g_{n-m+1}  & g_{n-m+2} & \cdots & g_{n-m-1} & g_{n-m}
\end{pmatrix}_{m\times n},
\end{equation}
where all the subscripts $j$ of $g_j$ are modulo $n$.
We call $M(g)_{m\times n}$ a {\em generalized circulant matrix} 
 associated with the polynomial $g(X)=\sum_{j=0}^{n-1}g_jX^j$.
In the special case when $m=n$,
$M(g)_{n\times n}$ is just the usual circulant matrix mentioned above.

For any $n$-th root of unity $\zeta=\omega^j$, we define a
column vector of length $m$:
\begin{equation}\label{GEV}
 V(\zeta)_m=(1,\zeta,\cdots,\zeta^{m-1})^T,
\end{equation}
where the superscript $T$ stands for the transpose.
A known elementary result says that $g(\omega^j)$
is an eigenvalue of the circulant matrix $M(g)_{n\times n}$
and $V(\omega^j)_n$ is a corresponding eigenvector,
i.e., $M(g)_{n\times n}V(\omega^j)_n=g(\omega^j)V(\omega^j)_n$.
An easy modification is as follows.

\begin{lemma}\label{lemma GEV}
For any $n$-th root $\omega^j$ of unity
and any polynomial $g(X)=\sum_{i=0}^{n-1}g_iX^i$ over $F$ with
degree $\deg g(X)<n$,
$$M(g)_{m\times n} V(\omega^j)_{n}=g(\omega^j) V(\omega^j)_m.$$
\end{lemma}

\pf Set $\zeta=\omega^j$. The $0$-th row of $M(g)_{m\times n} V(\omega^j)_{n}$,
which is a $m\times 1$ matrix, is the following element
$$
(g_0,g_1,\cdots,g_{n-1})
\begin{pmatrix}1\\ \zeta\\ \vdots \\ \zeta^{n-1} \end{pmatrix}
=g_0+g_1\zeta+\cdot+g_{n-1}\zeta^{n-1}=g(\zeta).
$$
Note that $\zeta^n=1$, hence $\zeta^{n-1}=\zeta^{-1}$.
We can get the $1$-th row of $M(g)_{m\times n} V(\omega^j)_{n}$ as follows:
$$
g_{n-1}+g_0\zeta+\cdot+g_{n-2}\zeta^{n-1}=
\zeta(g_0+g_1\zeta+\cdot+g_{n-1}\zeta^{-1})
=\zeta g(\zeta).
$$
And the $2$-th row is
$$
g_{n-2}+g_{n-1}\zeta+\cdot+g_{n-3}\zeta^{n-1}=
\zeta(g_{n-1}+g_0\zeta+\cdot+g_{n-2}\zeta^{n-1})
=\zeta^2 g(\zeta).
$$
Iterating it in this way, we get the equality of the lemma.
\qed

\medskip
Let $\Phi(\omega)_{m\times n}$ be the $m\times n$ matrix whose $j$-th column
for $j=0,1,\cdots,n-1$ is the column vector $V(\omega^j)_m$, i.e.
\begin{equation}\label{FTM}
\Phi(\omega)_{m\times n}
\!=\!\big( V(1)_m,\!V(\omega)_m,\cdots\!,V(\omega^{n-1})_m\!\big)
\!=\!\begin{pmatrix}
1 & 1 & \cdots  & 1 \\
1 & \omega &\cdots &  \omega^{n-1} \\
 \vdots & \vdots & \cdots &  \vdots\\
 1  &\!\omega^{m-1}\! & \cdots & \!\omega^{(n-1)(m-1)}\!
\end{pmatrix}.
\end{equation}
In the special case when $m=n$, $\Phi(\omega)_{n\times n}$
is just the known Fourier transform matrix for cyclic group of order $n$.
For convenience, we call $\Phi(\omega)_{m\times n}$ a {\em generalized
Fourier matrix}.

\begin{remark}\label{rank of Phi}\rm
We exhibit a typical argument in linear algebra.
Let $r=\min\{m,n\}$.
For any positive integer $k$ such that $k\leq r$,
any $k\times k$ sub-matrix taken from the first $k$ rows of
the matrix $\Phi(\omega)_{m\times n}$
is a Vandemond matrix formed by $k$ different elements,
hence is a non-degenerate matrix.
Thus, the rank of the matrix $\Phi(\omega)_{m\times n}$,
denoted by ${\rm rank}\,\Phi(\omega)_{m\times n}$,
is equal to $r=\min\{m,n\}$,
and the first $r$ rows of the matrix $\Phi(\omega)_{m\times n}$
are linearly independent.
\end{remark}

From Lemma \ref{lemma GEV}, we get an important formula.
By ${\rm diag}(a_0, a_1,\cdots,a_{n-1})$ we
denote the diagonal matrix with diagonal elements
$a_0,a_1,\cdots,a_{n-1}$.

\begin{lemma}\label{C-F}
For any polynomial $g(X)=\sum_{i=0}^{n-1}g_iX^i$ 
over $F$ with ${\deg g(X)<n}$,
\begin{equation}\label{C-F-e}
M(g)_{m\times n}\cdot\Phi(\omega)_{n\times n}
 =\Phi(\omega)_{m\times n}\cdot{\rm diag}
 \big(g(1),\,g(\omega),\,\cdots,\,g(\omega^{n-1})\big).
\end{equation}
\end{lemma}

\pf By the block multiplication of blocked matrices,
\begin{align*}
&M(g)_{m\times n}\cdot\Phi(\omega)_{n\times n}
=M(g)_{m\times n}\cdot\Big( V(1)_n,V(\omega)_n,\cdots,V(\omega^{n-1})_n\Big)\\
&=\Big(M(g)_{m\times n}V(1)_n,\,M(g)_{m\times n}V(\omega)_n,\,
\cdots,\,M(g)_{m\times n}V(\omega^{n-1})_n\Big)\\
&=\Big( g(1)V(1)_{m},\;g(\omega)V(\omega)_m,\;
 \cdots,\; g(\omega^{n-1})V(\omega^{n-1})_m \Big)\\
&=\Phi(\omega)_{m\times n}\cdot{\rm diag}
 \big(g(1),\,g(\omega),\,\cdots,\,g(\omega^{n-1})\big). \hskip30mm\qed
\end{align*}

As usual, $\gcd\big( f(X), g(X)\big)$ 
denotes the greatest common divisor.

\begin{theorem}\label{thm for circulant}
Let $g(X)$ be a polynomial over $F$ with $\deg g(X)<n$,  let
$d=\deg\gcd\big( g(X),\,X^n-1\big)$ and $r=\min\{m, n-d\}$.
Then ${\rm rank}\,M(g)_{m\times n}=r$ and
the first $r$ rows of $M(g)_{m\times n}$ are linearly independent.
\end{theorem}

\pf 
There are exactly $d$ indexes $j$ with $0\le j< n$ 
such that ${g(\omega^j)=0}$. 
So the matrix $\Phi(\omega)_{m\times n}\!\cdot\!
{\rm diag}\big(g(1),g(\omega),\cdots,g(\omega^{n-1})\big)$ of the
right hand side of Eqn \eqref{C-F-e} has exactly $n-d$ non-zero columns.
Assume that $j_1,\cdots,j_{n-d}$ 
are the indexes such that 
$g(\omega^{j_t})\ne 0$, $t=1,\cdots,n-d$.
Then the rank of $M(g)_{m\times n}\cdot\Phi(\omega)_{n\times n}$
is equal to the rank of the following $m\times(n-d)$ matrix
\begin{equation}\label{part}
\begin{array}{l}
\big(V(\omega^{j_1})_{m},\cdots, V(\omega^{j_{n-d}})_{m}\big)
\cdot {\rm diag}\big(g(\omega^{j_1}),\cdots,g(\omega^{j_{n-d}})\big)\\[5pt]
= \big(g(\omega^{j_1})V(\omega^{j_1})_{m},\;\cdots,
 g(\omega^{j_{n-d}})V(\omega^{j_{n-d}})_{m}\big).
\end{array}
\end{equation}
Similar to the argument in Remark \ref{rank of Phi}, we get that
$$
{\rm rank}\big(V(\omega^{j_1})_{m},\,\cdots,
 V(\omega^{j_{n-d}})_{m}\big) = r,
$$
and the first $r$ rows of the matrix
$\big(V(\omega^{j_1})_{m},\cdots, V(\omega^{j_{n-d}})_{m}\big)$ 
are linearly independent. Since 
${\rm diag}\big(g(\omega^{j_1}), \cdots, g(\omega^{j_{n-d}})\big)$ 
is non-degenerate,
the rank of the matrix in Eqn \eqref{part} is equal to $r$ and
the first $r$ rows of the matrix are linearly independent.
By Lemma \ref{C-F}, the rank of the matrix 
$M(g)_{m\times n}\!\cdot\!\Phi(\omega)_{n\times n}$
is equal to $r$ and the first $r$ rows of
$M(g)_{m\times n}\!\cdot\!\Phi(\omega)_{n\times n}$ are linearly independent.
Finally, since $\Phi(\omega)_{n\times n}$ is non-degenerate,
the rank of the matrix $M(g)_{m\times n}$
is equal to $r$ and the first $r$ rows of
$M(g)_{m\times n}$ are linearly independent.
\qed

\begin{corollary} Any consecutive $r$ rows of
$M(g)_{m\times n}$ are linearly independent,
where $r=\min\big\{m, n-\deg\gcd\big(g(X),X^n-1\big)\big\}$
is the rank of $M(g)_{m\times n}$.
\end{corollary}

\pf If $r=m$, the conclusion of the corollary is true trivially.
Assume $r=n-\deg\gcd\big(g(X),X^n-1\big)$.
Let $h$ be a multiple of $n$ such that $h\ge m$.
By Theorem \ref{thm for circulant},
${\rm rank}\, M(g)_{h\times n}=r$
and the first $r$ rows of $M(g)_{h\times n}$ are linearly independent.
When we shift circularly the rows of $M(g)_{h\times n}$,
the resultant matrix is still a generalized circulant matrix with rank $r$,
whose first $r$ rows are linearly independent.
Thus, any consecutive $r$ rows of $M(g)_{h\times n}$ are linearly independent.
Finally, any consecutive $r$ rows of $M(g)_{m\times n}$ are
consecutive $r$ rows from $M(g)_{h\times n}$, hence
are linearly independent.
\qed

\begin{remark}\label{why circulant}\rm
In the special case when $m=n$, as $\min\{n-d, m\}=n-d$,
we see that the rank of the circulant matrix $M(g)_{n\times n}$
is equal to $n-d$ and any consecutive $r$ rows of $M(g)_{n\times n}$
are linearly independent.
The former conclusion is well-known, e.g., \cite{AW}. Whereas,
the latter conclusion is important for cyclic codes
to construct generator matrices. Since this is a very original idea
to motivate our double circulant matrices,
we sketch it briefly as follows.

Let $F$ be a finite field such that $\gcd({\rm char} F, n)=1$.
In coding-theoretic notations, any $(a_0,a_1,\cdots,a_{n-1})\in F^n$
is called a {\em word} over $F$; any subspace $C$ of $F^n$ is called a
{\em linear code} over $F$, and the words in $C$ are called {\em code words}.
If $(g_{00},g_{01},\cdots,g_{0,n-1})$, $(g_{10},g_{11},\cdots,g_{1,n-1})$,
$\cdots$, $(g_{r-1,0},g_{r-1,1},\cdots,g_{r-1,n-1})$ are a basis of the
linear code $C$, then the matrix
$$
G=\begin{pmatrix}
 g_{00}& g_{01}& \cdots & g_{0,n-1}\\
 g_{10}& g_{11}& \cdots & g_{1,n-1}\\
 \cdots & \cdots & \cdots & \cdots \\
 g_{r-1,0}& g_{r-1,1}& \cdots & g_{r-1,n-1}\\
\end{pmatrix}_{r\times n}
$$
is said to be a {\em generator matrix} of the linear code $C$.
The generator matrices are useful for encoding and decoding.

Next, we consider the quotient ring $F[X]/\langle X^n-1\rangle$
of the polynomial ring $F[X]$ over the ideal $\langle X^n-1\rangle$
generated by $X^n-1$. Any polynomial
$a(X)=\sum_{j=0}^{n-1}a_jX^j\in F[X]/\langle X^n-1\rangle$
is identified with a word $(a_0,a_1,\cdots,a_{n-1})\in F^n$. In this way,
any ideal $C$ of $F[X]/\langle X^n-1\rangle$ is viewed
as a linear code, called a {\em cyclic code} over $F$.
Then any polynomial $g(X)=g_0+g_1X+\cdots+g_{n-1}X^{n-1}$ generates
a cyclic code $C=\big\{f(X)g(X)\pmod{X^n-1}\,\big|\,f(X)\in F[X]\big\}$
(and any cyclic code can be constructed in this way).
Since
$$\begin{array}{rcl}
 X\,g(X)\hskip-7pt&\equiv&\hskip-7pt
   g_{n-1}+g_0X+g_1X^2+\cdots+g_{n-2}X^{n-1}\pmod{X^n-1},\\
 X^2 g(X)\hskip-7pt&\equiv&\hskip-7pt
   g_{n-2}+g_{n-1}X+g_0X^2+\cdots+g_{n-3}X^{n-1}\pmod{X^n-1},\\
 \cdots~~~&\cdots&\qquad\cdots\qquad \cdots\qquad\cdots\qquad\cdots
\end{array}$$
we see that the row vectors of the following circulant matrix associated with $g(X)$
generate the linear code $C$:
$$
 M(g)_{n\times n}=\begin{pmatrix}g_0 & g_1 & \cdots &  g_{n-1} \\
 g_{n-1} & g_0 &\cdots &  g_{n-2} \\
 \cdots & \cdots & \cdots  & \cdots\\
 g_{1}  & g_{2} & \cdots  & g_0
\end{pmatrix}_{n\times n}.
$$
However, $M(g)_{n\times n}$ is not a generator matrix of $C$ in general
since the row vectors are linearly dependent in general.
Let $r=n-\deg\gcd\big(g(X),X^n-1\big)$.
By Theorem \ref{thm for circulant}, the first $r$ rows of $M(g)_{n\times n}$
form a basis of the cyclic code $C$. Taking the first $r$ rows, we get a
generator matrix of $C$ as follows:
$$
 G=\begin{pmatrix}g_0 & g_1 & \cdots &  g_{n-1} \\
 g_{n-1} & g_0 &\cdots &  g_{n-2} \\
 \cdots & \cdots & \cdots  & \cdots\\
 g_{n-r+1}  & g_{n-r+2} & \cdots  & g_{n-r}
\end{pmatrix}_{r\times n}.
$$
In particular, if $g(X)\,\big|\, (X^n-1)$, then $r=n-d$ where $d=\deg g(X)$;
and, writing $g(X)=g_0+g_1X+\cdots+g_{d}X^{d}$, we obtain a generator matrix
of $C$ as follows:
$$
 G=\begin{pmatrix}g_0 & g_1 & \cdots &  g_{d} \\
   & g_0 & g_1 &\cdots &  g_{d} \\
  &  & \cdots  & \cdots & \cdots & \cdots\\
  & & & g_0  & g_1 & \cdots  & g_{d}
\end{pmatrix}_{r\times n}.
$$

\end{remark}

\section{Double circulant matrices}\label{d-c matrix}

In this section we always assume that:

\hangindent10mm
$\bullet$~ $m$, $n$ and $n'$ are positive integers, 
and $\ell=\gcd(n,n')$;

\hangindent10mm
$\bullet$~  $F$ is a field with
${\rm char} F$ being zero or coprime to both $n$ and $n'$;

\hangindent10mm
$\bullet$~ $\omega$ and $\omega'$, respectively,
are a primitive $n$-th and a primitive $n'$-th, respectively,
roots of unity in a suitable extension of $F$.

Let $g(X)=\sum_{j=0}^{n-1}g_jX^j$ and $g'(X)=\sum_{j'=0}^{n'-1}g_{j'}X^{j'}$
be polynomials over $F$ with $\deg g(X)<n$ and $\deg g'(X)<n'$.
For any positive integer $m$, by Eqn \eqref{GCM}, we have
the generalized circulant matrices $M(g)_{m\times n}$ associated with $g(X)$,
and the generalized circulant matrices $M(g')_{m\times n'}$ associated with $g'(X)$.

\begin{definition}\label{def d-circulant}\rm
Concatenating  $M(g)_{m\times n}$ and $M(g')_{m\times n'}$ side-by-side,
we get a $m\times(n+n')$ matrix
\begin{equation}\label{GDC}
M(g,g')_{m\times(n+n')}=\Big(M(g)_{m\times n}\,\big|\,M(g')_{m\times n'}\Big),
\end{equation}
which we call by a {\em double circulant matrix}
associated with the polynomials $g(X)$ and $g'(X)$.
\end{definition}

Note that, in the special case when $m=n+n'$, we have
a square double circulant matrix
$M(g,g')_{(n+n')\times(n+n')}$;
however, we'll see that this square matrix is not full-rank
for any $g(X)$ and $g'(X)$,
see Corollary \ref{not full-rank} below.

\begin{example}\rm
Let $F$ be the complex field,
$n=2$, $n'=3$, $g(X)=-1+X$, $g'(X)=-2+X+X^2$. Then for $m=6$ we have
$$
M(g,g')_{6\times(2+3)}=
\begin{pmatrix}
  -1&1&-2&1&1\\ 1&-1&1&-2&1\\-1&1&1&1&-2\\
  1&-1&-2&1&1\\ -1&1&1&-2&1\\1&-1&1&1&-2
\end{pmatrix}.
$$
While for $m=5$, we see that $M(g,g')_{5\times(2+3)}$ is a square matrix
of size $5$ formed by the first $5$ rows of $M(g,g')_{6\times(2+3)}$.
One can check that both $M(g,g')_{5\times(2+3)}$
and $M(g,g')_{6\times(2+3)}$ have rank $3$, and the first $3$ rows
of $M(g,g')_{5\times(2+3)}$ are linearly independent.
\end{example}

\begin{example}\rm
Let $F$ be the complex field,
$n=4$, $n'=2$, $g(X)=-2+X+X^2$, $g'(X)=-1+X$. Then for $m=4$ we have
$$
M(g,g')_{4\times(4+2)}=
\begin{pmatrix}
  -2&1&1&0&-1&1\\ 0&-2&1&1&1&-1\\1&0&-2&1&-1&1&\\ 1&1&0&-2&1&-1
\end{pmatrix}.
$$
While for $m=6$, we see that $M(g,g')_{6\times(4+2)}$ is a square matrix
of size $6$, which has the first $4$ rows the same as $M(g,g')_{4\times(4+2)}$,
and has the last $2$ rows being a copy of the first $2$ rows of $M(g,g')_{4\times(4+2)}$.
One can check that both $M(g,g')_{4\times(4+2)}$
and $M(g,g')_{6\times(4+2)}$ have rank $3$, and the first $3$ rows
of $M(g,g')_{4\times(4+2)}$ are linearly independent.
\end{example}

Our aim is to find the rank of any double circulant matrix and
find a maximal linearly independent set of rows of the matrix.

Let $V(\zeta)_m$ be defined as in Eqn\eqref{GEV}.
Then both $\begin{pmatrix}V(\zeta)_{n}\\ 0_{n'\times 1} \end{pmatrix}$
and $\begin{pmatrix}0_{n\times 1}\\ V(\zeta)_{n'} \end{pmatrix}$
are $(n+n')\times 1$ matrces,
where $0_{n\times 1}$ denotes the $n\times 1$ zero matrix.

\begin{lemma}\label{lemma GEV-d}
For any $n$-th root $\omega^j$ of unity
and $n'$-th root $\omega'^{j'}$ of unity,
any polynomial $g(X)=\sum_{k=0}^{n-1}g_kX^k$ and
$g'(X)=\sum_{k'=0}^{n'-1}g_{k'}X^{k'}$ over $F$ with
degree $\deg g(X)<n$ and $\deg g'(X)<n'$,
\begin{equation}\label{GEV-d}
M(g,g')_{m\times(n+n')}
\begin{pmatrix}V(\omega^j)_{n}\\ 0_{n'\times 1} \end{pmatrix}
 = g(\omega^j) V(\omega^j)_m;
\end{equation}
\begin{equation}\label{GEV-d'}
M(g,g')_{m\times(n+n')}
\begin{pmatrix}0_{n\times 1}\\V(\omega'^{j'})_{n'} \end{pmatrix}
 = g'(\omega'^{j'}) V(\omega'^{j'})_m.
\end{equation}
\end{lemma}

\pf
By the block multiplication of blocked matrices, we have:
\begin{align*}
&M(g,g')_{m\times(n+n')}
\begin{pmatrix}V(\omega^j)_{n}\\ 0_{n'\times 1}\end{pmatrix}
=\Big(M(g)_{m\times n}\,\big|\,M(g')_{m\times n'}\Big)
\begin{pmatrix}V(\omega^j)_{n}\\ 0_{n'\times 1}\end{pmatrix}\\
&=M(g)_{m\times n}V(\omega^j)_n+M(g')_{m\times n'}\cdot 0_{n'\times 1}
=M(g)_{m\times n}V(\omega^j)_n\\
&=g(\omega^j) V(\omega^j)_m,
\end{align*}
where the last equality is obtained by Lemma \ref{lemma GEV}. And:
\begin{align*}
&M(g,g')_{m\times(n+n')}
\begin{pmatrix} 0_{n\times 1}\\ V(\omega'^{j'})_{n'}\end{pmatrix}
=\Big(M(g)_{m\times n}\,\big|\,M(g')_{m\times n'}\Big)
\begin{pmatrix} 0_{n\times 1}\\ V(\omega'^{j'})_{n'}\end{pmatrix}\\
&=M(g)_{m\times n}\cdot 0_{n\times 1}+M(g')_{m\times n'}V(\omega'^{j'})_{n'}
=M(g')_{m\times n'}V(\omega'^{j'})_{n'}\\
&=g'(\omega'^{j'}) V(\omega'^{j'})_m,
\end{align*}
where the last equality is still by Lemma \ref{lemma GEV}.
\qed


As in Eqn \eqref{FTM}, we have the generalized Fourier matrices
$\Phi(\omega)_{m\times n}$ and $\Phi(\omega')_{m\times n'}$.
Similarly to Eqn \eqref{GDC}, concatenating them side-by-side,
we get a $m\times(n+n')$ matrix:
\begin{equation}\label{FTM-d}
\Phi(\omega,\omega')_{m\times(n+n')}=
\Big(\Phi(\omega)_{m\times n}\,\big|\,\Phi(\omega')_{m\times n'}\Big).
\end{equation}
On the other hand, with the square matrices
$\Phi(\omega)_{n\times n}$ and $\Phi(\omega')_{n'\times n'}$,
we construct a diagonal blocked matrix
\begin{equation}\label{GFTM}
{\rm diag}\Big(\Phi(\omega)_{n\times n},\,\Phi(\omega')_{n'\times n'}\Big)
=\begin{pmatrix}\Phi(\omega)_{n\times n}\\[5pt]
 &\Phi(\omega')_{n'\times n'} \end{pmatrix}.
\end{equation}

Similarly to Lemma \ref{C-F}, from Lemma \ref{lemma GEV-d}
we get the following important formula for double circulant matrices.

\begin{lemma}\label{C-F-d}
For any polynomial $g(X)=\sum_{k=0}^{n-1}g_kX^k$ and
$g'(X)\!=\!\sum_{k'=0}^{n'-1}g_{k'}X^{k'}$ over $F$ with
degree $\deg g(X)<n$ and $\deg g'(X)<n'$,
\begin{align*}
&M(g,g')_{m\times(n+n')}\cdot
 {\rm diag}\Big(\Phi(\omega)_{n\times n},\Phi(\omega')_{n'\times n'}\Big)=\\
&~ \Phi(\omega,\omega')_{m\times(n+n')}\cdot
 {\rm diag}\big(g(1),g(\omega),\cdots,g(\omega^{n-1}),
   g'(1),g'(\omega'),\cdots,g'(\omega'^{n'-1})\big).
\end{align*}
\end{lemma}

\pf To shorten the notations, we denote
\begin{equation*}
\begin{array}{l}
 D(\omega,\omega')={\rm diag}\Big(\Phi(\omega)_{n\times n},\,\Phi(\omega')_{n'\times n'}\Big),\\
 D(g,g')={\rm diag}\big(g(1),g(\omega),\cdots,g(\omega^{n-1}),\,
   g'(1),g'(\omega'),\cdots,g'(\omega'^{n'-1})\big).
\end{array}
\end{equation*}
Then the formula we are proving is
\begin{equation}\label{MD=}
M(g,g')_{m\times(n+n')}\cdot D(\omega,\omega')
 =\Phi(\omega,\omega')_{m\times(n+n')}\cdot D(g,g').
\end{equation}
Take any column $L$ of $D(\omega,\omega')$ and multiply it with
$M(g,g')_{m\times(n+n')}$. There are two cases:

If the column $L$ is located within $\Phi(\omega)_{n\times n}$, then it is the form
$L=\begin{pmatrix}V(\omega^j)_{n}\\ 0_{n'\times 1} \end{pmatrix}$;
by Eqn \eqref{GEV-d}, $M(g,g')_{m\times(n+n')}\cdot L=g(\omega^j)V(\omega^j)_m$,
which is just the corresponding column of
$\Phi(\omega,\omega')_{m\times(n+n')}\cdot D(g,g')$.

Otherwise, $L$ is located within $\Phi(\omega')_{n'\times n'}$ and is the form
$L=\begin{pmatrix}0_{n\times 1}\\ V(\omega'^{j'})_{n'} \end{pmatrix}$;
by Eqn \eqref{GEV-d'},
$M(g,g')_{m\times(n+n')}\cdot L=g'(\omega'^{j'})V(\omega'^{j'})_m$,
which is the corresponding column of $\Phi(\omega,\omega')_{m\times(n+n')}\cdot D(g,g')$.
\qed

\begin{theorem}\label{thm d-circulant}
Let $M(g,g')_{m\times(n+n')}$ be the double circulant matrix
associated with polynomials $g(X)=\sum_{k=0}^{n-1}g_kX^k$ and
$g'(X)=\sum_{k'=0}^{n'-1}g_{k'}X^{k'}$ over $F$ of
degree $\deg g(X)<n$ and $\deg g'(X)<n'$.
Let $\ell=\gcd(n,n')$,
\begin{equation}\label{d=}\textstyle
 d=\deg\frac{\gcd\!\big(g(X),\,X^n-1\big)
 \cdot\gcd\!\big(g'(X),\,X^{n'}-1\big)\cdot\big(X^\ell -1\big)}
     {\gcd\!\big(g(X)g'(X),\,X^\ell-1\big)},
\end{equation}
and $r=\min\{m,\,n+n' - d\}$. Then
$${\rm rank}\,M(g,g')_{m\times(n+n')}=r, $$
and the first $r$ rows of $M(g,g')_{m\times(n+n')}$ are linearly independent.
\end{theorem}

\pf
Since $D(\omega,\omega')
={\rm diag}\big(\Phi(\omega)_{n\times n},\Phi(\omega')_{n'\times n'}\big)$
is non-degenerate, by the formula in Lemma \ref{C-F-d}
(see Eqn \eqref{MD=} for the shortened notation),
we just need to show that

$\bullet$~ {\it The rank of the matrix
$\Phi(\omega,\omega')_{m\times(n+n')}\cdot D(g,g')$ is equal to $r$, and
the first $r$ rows of the matrix are linearly independent.}

An element $\zeta$ in any extension of $F$
is a root of both $X^n-1$ and $X^{n'}-1$ if and only if $\zeta$ is a root
of $X^\ell-1$. Thus,
\begin{equation}\label{gcd n n'}
 \gcd\big(X^n-1,\,X^{n'}-1\big)=X^{\ell}-1.
\end{equation}
For convenience, in the following we denote
\begin{equation}\label{e e'}
\begin{array}{l}
 e=\deg\gcd\!\big( g(X),X^n-1\big),~~
 e'=\deg\gcd\!\big( g'(X),X^{n'}-1\big),\\[3pt]
\displaystyle
 \bar e=\deg\frac{X^\ell -1}{\gcd\!\big( g(X)g'(X),~X^\ell-1\big)}.
\end{array}
\end{equation}
By Eqn \eqref{d=}, $d=e+e'+\bar e$.
Note that $\eta=\omega^{\frac{n}{\ell}}$ is a primitive $\ell$-th root of unity,
hence $\eta^0=1,\eta,\cdots,\eta^{\ell-1}$ are all roots of $X^\ell-1$.
And note that $\eta^0=1,\eta,\cdots,\eta^{\ell-1}$ are also roots of $X^{n'}-1$,
cf. Eqn \eqref{gcd n n'}.
So there are exactly $\bar e$ roots of $X^\ell -1$,
say $\eta^{\bar j_1},\cdots,\eta^{\bar j_{\bar e}}$ where
$0\le\bar j_1<\cdots<\bar j_{\bar e}<\ell$, such that 
\begin{equation}\label{bar j}
 g(\eta^{\bar j_i})\ne 0\ne g'(\eta^{\bar j_i}),~~~~~ i=1,\cdots,\bar e.
\end{equation}

The key step of the proof is to determine which of the columns of
the matrix $\Phi(\omega,\omega')_{m\times(n+n')}\cdot D(g,g')$
contribute essentially to the rank of the matrix.
They are determined by two observations.

First, if $g(\omega^{j})=0$ (or $g'(\omega'^{j'})=0$), then
the corresponding column of the matrix
$\Phi(\omega,\omega')_{m\times(n+n')}\cdot D(g,g')$
is $g(\omega^j)V(\omega^j)_m=0$, hence should be ignored.
So we have indexes $0\le j_1<\cdots<j_{n-e-\bar e}<n$ and
$0\le j'_1<\cdots<j'_{n'-e'-\bar e}<n'$ such that
\begin{equation}\label{non-zeros}
\begin{matrix}
 g(\omega^{j_1}),~\cdots,~g(\omega^{j_{n-e-\bar e}}),~
 g(\eta^{\bar j_1}),~\cdots,~g(\eta^{\bar j_{\bar e}});\\[5pt]
 g'(\omega'^{j'_1}),~\cdots,~g'(\omega'^{j'_{n'-e'-\bar e}}),~
 g'(\eta^{\bar j_1}),~\cdots,~g'(\eta^{\bar j_{\bar e}});
\end{matrix}
\end{equation}
are all non-zero and the corresponding columns of
$\Phi(\omega,\omega')_{m\times(n+n')}\cdot D(g,g')$
are considered to find the rank of the matrix.

Next, the column vector $V(\eta^{\bar j_i})_m$
of $\Phi(\omega,\omega')_{m\times(n+n')}$ for $i=1,\cdots,\bar e$,
which we selected to compute the rank of the matrix,
appears in the list \eqref{non-zeros} twice.
However, it contributes only one to the rank of
$\Phi(\omega,\omega')_{m\times(n+n')}\cdot D(g,g')$.
So, for each $i=1,\cdots,\bar e$ we can take only one column
$V(\eta^{\bar j_i})_m$ from $\Phi(\omega,\omega')_{m\times(n+n')}$
to compute the rank of
$\Phi(\omega,\omega')_{m\times(n+n')}\cdot D(g,g')$.

According to the above two observations, we rearrange columns of
the matrix $\Phi(\omega,\omega')_{m\times(n+n')}\cdot D(g,g')$
and construct a sub-matrix of it with the selected columns as follows:
$$\begin{array}{l}
\Phi=
\Big(V(\omega^{j_1})_m,\cdots,V(\omega^{j_{n-e-\bar e}})_m,\\
 \hskip11mm
 V(\eta^{\bar j_1})_m,\cdots,V(\eta^{\bar j_{\bar e}})_m,\,
 V(\omega'^{j'_1})_m,\cdots,V(\omega'^{j'_{n'-e'-\bar e}})_m \Big),\\[5pt]
D={\rm diag}
\big( g(\omega^{j_1}),\cdots, g(\omega^{j_{n-e-\bar e}}),\\
 \hskip17mm
 g(\eta^{\bar j_1}),\cdots, g(\eta^{\bar j_{\bar e}}),\,
 g'(\omega'^{j'_1}),\cdots,g'(\omega'^{j'_{n'-e'-\bar e}})\big).
\end{array}
$$
Then
\begin{equation}\label{Phi D}
{\rm rank}\big(\Phi\!\cdot\! D\big)=
{\rm rank}\big(\Phi(\omega,\omega')_{m\times(n+n')}\cdot D(g,g')\big);
\end{equation}
and, by Eqn \eqref{gcd n n'}, the following elements
are distinct from each other:
$$
\omega^{j_1},\cdots, \omega^{j_{n-e-\bar e}},\,
\eta^{\bar j_1},\cdots, \eta^{\bar j_{\bar e}},\,
\omega'^{j'_1},\cdots, \omega'^{j'_{n'-e'-\bar e}}.
$$

Finally, $D$ is non-degenerate,
and $\Phi$ is an $m\times(n+n'-d)$ matrix such that,
for any positive integer $h$ with $h\le r=\min\{m,\,n+n'-d\}$,
any $h\times h$ sub-matrix taken from the first $h$ rows of $\Phi$
is a Vandemond matrix formed by $h$ different elements.
By the same argument in Remark \ref{rank of Phi}, we get that
${\rm rank}\big(\Phi\!\cdot\! D\big)=r$, and the first $r$ rows of
the matrix $\Phi\!\cdot\! D$ are linearly independent. So
${\rm rank}\big(\Phi(\omega,\omega')_{m\times(n+n')}\!\cdot\! D(g,g')\big)=r$
and the first $r$ rows of the matrix
$\Phi(\omega,\omega')_{m\times(n+n')}\cdot D(g,g')$
are linearly independent.
\qed

\begin{corollary}\label{cor thm for d-circulant}
Notations are the same as in Theorem \ref{thm d-circulant}.
Any consecutive $r$ rows of
$M(g,g')_{m\times(n+n')}$ are linearly independent.
\end{corollary}

\pf 
If $r=m$, there is nothing to do.
Let $h$ be a common multiple of $n$ and $n'$ such that $h\ge m$.
By Theorem \ref{thm d-circulant},
$${\rm rank}\, M(g,g')_{h\times(n+n')}=r, $$
and the first $r$ rows of $M(g,g')_{h\times(n+n')}$ are linearly independent.
When we shift circularly the rows of the matrix $M(g,g')_{h\times(n+n')}$,
the resultant matrix is still a double circulant matrix with rank $r$,
whose first $r$ rows are linearly independent.
Thus, any consecutive $r$ rows of $M(g,g')_{h\times(n+n')}$ are linearly independent.
Any consecutive $r$ rows of $M(g,g')_{m\times(n+n')}$ are
consecutive $r$ rows from $M(g,g')_{h\times(n+n')}$, hence
are linearly independent.
\qed

\begin{corollary}\label{not full-rank}
${\rm rank}\,M(g,g')_{m\times(n+n')}<n+n'$.
\end{corollary}

\pf Since $\ell= \gcd(n,n')\ge 1$, by Eqn \eqref{d=} it is easy to check
 that $d\ge\ell\ge 1$. \qed

\smallskip
If we consider the square double circulant matrix
$M(g,g')_{(n+n')\times(n+n')}$ where we take $m=n+n'$,
then we can sow a little more information.

\begin{corollary}\label{eigen 0}
Let notations be as in Theorem \ref{thm d-circulant}.
Assume that $m=n+n'$.
Then $0$ is an eigenvalue of the matrix $M(g,g')_{m\times m}$
and the following $d$ vectors are eigenvectors of the eigenvalue $0$
which are linearly independent:

{\rm(E1)}~ the vectors
$\begin{pmatrix}V(\omega^j)_n\\ 0_{n'\times 1}\end{pmatrix}$
for every root $\omega^j$ of $\gcd\big(g(X), X^n-1\big)$;

{\rm(E2)}~ the vectors
$\begin{pmatrix}0_{n\times 1}\\ V(\omega'^{j'})_{n'}\end{pmatrix}$
for every root $\omega'^{j'}$ of $\gcd\big(g'(X), X^{n'}-1\big)$;

\hangindent15mm
{\rm(E3)}~ the vectors
$\begin{pmatrix}-g'(\eta^{\bar j})V(\eta^{\bar j})_n\\
   g(\eta^{\bar j})V(\eta^{\bar j})_{n'}\end{pmatrix}$
where $\eta=\omega^{\frac{n}{\ell}}$ and $0\le {\bar j}<\ell$
such that $g(\eta^{\bar j})\ne 0\ne g'(\eta^{\bar j})$.
\end{corollary}

\pf Let $e$, $e'$ and $\bar e$ be as in Eqn \eqref{e e'},
hence $d=e+e'+\bar e$.
If $g(\omega^j)=0$, by Lemma \ref{lemma GEV-d},
$$
 M(g,g')_{m\times m}
 \begin{pmatrix}V(\omega^j)_n\\ 0_{n'\times 1}\end{pmatrix}
 =0\cdot V(\omega^j)_m=0\cdot
 \begin{pmatrix}V(\omega^j)_n\\ 0_{n'\times 1}\end{pmatrix}.
$$
So the vectors in (E1) are eigenvectors of the eigenvalue $0$ of
$M(g,g')_{m\times m}$. The number of the vectors in (E1) is equal to $e$.
Similarly, there are $e'$ vectors in (E2) which are eigenvectors
of the eigenvalue $0$.
Finially, there are exactly $\bar e$ roots
$\eta^{\bar j}$ of $X^\ell -1$ for $\bar j=\bar j_i$,
with $i=1,\cdots,{\bar e}$, such that
$g(\eta^{\bar j})\ne 0\ne g'(\eta^{\bar j})$,
see Eqn \eqref{bar j}.
And by block multiplication of blocked matrices,
\begin{align*}
&M(g,g')_{m\times m}
\begin{pmatrix}-g'(\eta^{\bar j})V(\eta^{\bar j})_n\\
   g(\eta^{\bar j})V(\eta^{\bar j})_{n'}\end{pmatrix}
\!=\!\Big(M(g)_{m\times n}\,\big|\,M(g')_{m\times n'}\Big)
\begin{pmatrix}-g'(\eta^{\bar j})V(\eta^{\bar j})_n\\
   g(\eta^{\bar j})V(\eta^{\bar j})_{n'}\end{pmatrix}\\
&=-g'(\eta^{\bar j})M(g)_{m\times n}V(\eta^{\bar j})_n
 +g(\eta^{\bar j})M(g')_{m\times n'}V(\eta^{\bar j})_{n'}\\
&=-g'(\eta^{\bar j})g(\eta^{\bar j}) V(\eta^{\bar j})_m
 +g(\eta^{\bar j})g'(\eta^{\bar j}) V(\eta^{\bar j})_m =0.
\end{align*}
That is, the $\bar e$ vectors in (E3) are
eigenvectors of the eigenvalue $0$.
Finally, it is easy to check that
the above $d$ vectors are linearly independent.
\qed

\section{Applications to coding theory}\label{coding theory}

Continuing Remark \ref{why circulant}, we apply the
results on double circulant matrices
to two questions in coding theory. The materials also illustrate why
we are concerned with double circulant matrices.

\subsection{Quasi-cyclic codes of index $1\frac{1}{2}$}

Let $F$ be a finite field, and $n$ be a positive integer
with $\gcd(n,{\rm char}\, F)=1$.
We still consider the quotient ring $R=F[X]/\langle X^n-1\rangle$.
The product
$$F[X]/\langle X^n-1\rangle\times F[X]/\langle X^n-1\rangle$$
is a free $R$-module of rank $2$. In coding-theoretic notations,
any element of the module, $\big(a(X),a'(X)\big)$ with
$a(X)=\sum_{j=0}^{n-1}a_jX^j$ and $a'(X)=\sum_{j=0}^{n-1}a'_jX^j$,
is viewed as a word $(a_0,\cdots,a_{n-1},\,a'_0,\cdots,a'_{n-1})\in F^{2n}$,
and any $R$-submodule $C$ of the module
is said to be a {\em quasi-cyclic code of index $2$}.
More generally, any $R$-submodule of the
free $R$-module $R^m$ of rank $m$ is said to be a
{\em quasi-cyclic code of index $m$}.
So, cyclic codes are just quasi-cyclic codes of index $1$.

Though cyclic codes are a 
widely studied class of linear codes,
it is still a long-standing open question (cf. \cite{MW}):
whether or not the cyclic codes are asymptotically good?

However, it has been known since 1960's (see \cite{CPW}) that
the quasi-cyclic codes of index $2$ are asymptotically good.
In fact, for any integer $m>1$,
the quasi-cyclic codes of index $m$ are asymptotically good,
e.g., see \cite{FL}, \cite{LS}.

Thus, in \cite{FH} we introduced the quasi-cyclic codes of index
$1\frac{1}{2}$, and proved that they are asymptotically good.
We sketch it as follows.

Assume that $n$ is even. We consider the following
$F[X]/\langle X^n-1\rangle$-module: 
\begin{equation}\label{module}
F[X]/\langle X^n-1\rangle\times F[X]/\langle X^{\frac{n}{2}}-1\rangle.
\end{equation}
We name any submodule $C$ of the $F[X]/\langle X^n-1\rangle$-module
\eqref{module} by a {\em quasi-cyclic code of index $1\frac{1}{2}$}.
For any element of the module \eqref{module}:
$$ \big(g(X),g'(X)\big),\quad{\rm where}~
g(X)=\sum_{j=0}^{n-1}g_jX^j, ~ g'(X)=\sum_{j'=0}^{\frac{n}{2}-1}g'_{j'}X^{j'},
$$
a quasi-cyclic code $C_{g,g'}$ of index $1\frac{1}{2}$ can be constructed
as follows:
$$
C_{g,g'}\!=\!\big\{\big(f(X)g(X)\,({\rm mod}\,X^n-1),\,
  f(X)g'(X)\,({\rm mod}\,X^{\frac{n}{2}}-1)\big)\,\big|\, f(X)\in F[X] \big\},
$$
which we call a quasi-cyclic code of index $1\frac{1}{2}$ {\em generated by
$\big(g(X),\,g'(X)\big)$}.
Then we proved that this kind of
quasi-cyclic codes $C_{g,g'}$ of index $1\frac{1}{2}$
is asymptotically good.

A natural question is: how to get a generator matrix of the code $C_{g,g'}$?

\begin{table}[h]
\begin{center}
$\def\arraystretch{1.4}
 \begin{array}{c|c}
 \mbox{Elements of \eqref{module}} & \mbox{Corresponding code words}\\ \hline
 \big(g(X),g'(X)\big) & (g_0,g_1,\cdots,g_{n-1},\,g'_0,g'_1,\cdots,g'_{\frac{n}{2}-1})\\
 \big(Xg(X),Xg'(X)\big) &
  (g_{n-1},g_0,\cdots,g_{n-2},\,g'_{\frac{n}{2}-1},g'_0,\cdots,g'_{\frac{n}{2}-2})\\
  \cdots~~~~\cdots &\cdots~~~~~~~~\cdots\\
 \big(X^{\frac{n}{2}-1}g(X),X^{\frac{n}{2}-1}g'(X)\big) &
  (g_{\frac{n}{2}+1},g_{\frac{n}{2}+2},\cdots,g_{\frac{n}{2}},\,
  g'_{1},g'_2,\cdots,g'_{0})\\
 \big(X^{\frac{n}{2}}g(X),X^{\frac{n}{2}}g'(X)\big) &
  (g_{\frac{n}{2}},g_{\frac{n}{2}+1},\cdots,g_{\frac{n}{2}-1},\,
  g'_{0},g'_1,\cdots,g'_{\frac{n}{2}-1})\\
  \cdots~~~~\cdots &\cdots~~~~~~~~\cdots\\
 \big(X^{n-1}g(X),X^{n-1}g'(X)\big) &
  (g_{1},g_{2},\cdots,g_{0},\, g'_{1},g'_2,\cdots,g'_{0})\\
\end{array}$
\caption{Generating code words of $C_{g,g'}$}\label{element and word}
\end{center}
\end{table}

From Table 1, a double circulant matrix comes up:
$$
M(g,g')_{n\times(n+\frac{n}{2})}=
\begin{pmatrix}
g_0&g_1&\cdots&g_{n-1}&g'_0&g'_1&\cdots&g'_{\frac{n}{2}-1}\\
g_{n-1}&g_0&\cdots&g_{n-2}&g'_{\frac{n}{2}-1}&g'_0&\cdots&g'_{\frac{n}{2}-2}\\
\cdots & \cdots & \cdots & \cdots & \cdots & \cdots & \cdots & \cdots\\
g_{\frac{n}{2}+1}&g_{\frac{n}{2}+2}&\cdots&g_{\frac{n}{2}}&
  g'_{1} & g'_2 & \cdots & g'_{0}\\
g_{\frac{n}{2}}&g_{\frac{n}{2}+1}&\cdots&g_{\frac{n}{2}-1}&
  g'_{0} & g'_1 & \cdots & g'_{\frac{n}{2}-1}\\
\cdots & \cdots & \cdots & \cdots & \cdots & \cdots & \cdots & \cdots\\
g_{1} & g_{2} & \cdots & g_{0} & g'_{1} & g'_2 & \cdots & g'_{0}\end{pmatrix},
$$
whose rows generate the linear code $C_{g,g'}$.
However, $M(g,g')_{n\times(n+\frac{n}{2})}$
is not a generator matrix of $C_{g,g'}$ because its rows are linearly dependent
in general.

\begin{theorem}\label{q-cyclic 1.5}
Let $n$ be an even positive integer and $F$ be a finite field with
characteristic coprime to $n$.
Let $C_{g,g'}$ be the quasi-cyclic code of index $1\frac{1}{2}$ generated by
$\big(g(X),g'(X)\big)$ where $g(X)=\sum_{j=0}^{n-1}g_jX^j$
and $g'(X)=\sum_{j'=0}^{\frac{n}{2}-1}g'_{j'}X^{j'}$. Let
$$r=n-\deg\Big(\gcd\big(g(X),X^{\frac{n}{2}}+1\big)\cdot
 \gcd\big(g(X),g'(X),X^{\frac{n}{2}}-1\big)\Big).
$$
Then the dimension $\dim C_{g,g'}=r$ and the first $r$ rows of the double circulant
matrix $M(g,g')_{n\times(n+\frac{n}{2})}$ associated with $g(X)$ and $g'(X)$
form a generator matrix of $C_{g,g'}$.
\end{theorem}

\pf We begin with a remark on Eqn \eqref{d=}.
By a similar argument for Eqn \eqref{gcd n n'}
(i.e., checking roots of polynomials), for any $g(X)$ and $g'(X)$ we can show that
\begin{equation}\label{d==}
\begin{array}{l}
 \frac{\gcd\!\big(g(X),\,X^n-1\big)
 \cdot\gcd\!\big(g'(X),\,X^{n'}-1\big)}
     {\gcd\!\big(g(X)g'(X),\,X^\ell-1\big)}\\
~ =\gcd\!\big( g(X),\frac{X^n-1}{X^\ell-1}\big)
  \!\cdot\!\gcd\!\big( g'(X),\frac{X^{n'}-1}{X^\ell-1}\big)
  \!\cdot\!\gcd\!\big( g(X),g'(X),X^\ell-1\big).
\end{array}
\end{equation}

Return to the assumption of the theorem.
We have shown that the rows of the matrix $M(g,g')_{n\times(n+\frac{n}{2})}$
generate $C_{g,g'}$.
Apply Theorem \ref{thm d-circulant} and Eqn~\eqref{d==}
to $M(g,g')_{n\times(n+\frac{n}{2})}$.
Since $n'=\frac{n}{2}$ hence $\ell=\gcd(n,\frac{n}{2})=\frac{n}{2}=n'$,
we see that ${\rm rank}\,M(g,g')_{n\times(n+\frac{n}{2})}=r$
(cf. Eqn~\eqref{d=} and Eqn~\eqref{d==}),
and the first $r$ rows of $M(g,g')_{n\times(n+\frac{n}{2})}$ are
linearly independent. Thus $\dim C_{g,g'}=r$ and the first
$r$ rows of $M(g,g')_{n\times(n+\frac{n}{2})}$ form a basis of
$C_{g,g'}$.
\qed

\smallskip
The first conclusion ``$\dim C_{g,g'}=r$'' of the theorem has been
obtained in \cite{FH} but proved in another way.

\begin{example}\rm
Let $F={\Bbb Z}_3$ (the integer residue ring modulo $3$)
be the finite field with $3$ elements,
$n=4$, $n'=2$, $g(X)=1+X+X^2$, $g'(X)=2+X$.
Then $\gcd\big(g(X), X^2+1\big)=1$, $\gcd\big(g(X), g'(X), X^2-1\big)=1+X$.
By Theorem \ref{q-cyclic 1.5}, $\dim C_{g,g'}=3$.
From the double circulant matrix
$$
M(g,g')_{4\times(2+4)}=
\begin{pmatrix}
  1&1&1&0&2&1\\ 0&1&1&1&1&2\\1&0&1&1&2&1\\ 1&1&0&1&1&2
\end{pmatrix},
$$
we get a generator matrix $G$ of $C_{g,g'}$ as follows:
$$
G=
\begin{pmatrix}
  1&1&1&0&2&1\\ 0&1&1&1&1&2\\1&0&1&1&2&1
\end{pmatrix}.
$$
\end{example}

\subsection{Double cyclic codes}

Let $n$, $n'$ be positive integers, and $F$ be a finite field with characteristic
coprime to both $n$ and $n'$. Consider the following product
\begin{equation}\label{d-module}
F[X]/\langle X^n-1\rangle\times F[X]/\langle X^{n'}-1\rangle
\end{equation}
as an $F[X]$-module. Any $F[X]$-submodule $C$ of the
$F[X]$-module \eqref{d-module}
is said to be a {\em double cyclic code}, see \cite{ACT}.
In the general case, note that, \eqref{d-module} is neither an
$F[X]/\langle X^n-1\rangle$-module nor an
$F[X]/\langle X^{n'}-1\rangle$-module;
it can be viewed as an $F[X]/\langle X^{h}-1\rangle$-module where
$h$ is the least common multiple of $n$ and $n'$.

Any element of the $F[X]$-module \eqref{d-module}:
$$ \big(g(X),g'(X)\big),\quad{\rm where}~
g(X)=\sum_{j=0}^{n-1}g_jX^j, ~ g'(X)=\sum_{j'=0}^{n'-1}a'_{j'}X^{j'},
$$
generates a double cyclic code $C_{g,g'}$ as follows:
$$
C_{g,g'}=\big\{\big(f(X)g(X)\,({\rm mod}\,X^n-1),\,
  f(X)g'(X)\,({\rm mod}\,X^{n'}-1)\big)\,\big|\, f(X)\in F[X] \big\}.
$$
From Theorem \ref{thm d-circulant}, we can obtain straightforwardly
the dimension $r$ of $C_{g,g'}$ and a generator matrix $G$ of $C_{g,g'}$
as follows.

\begin{theorem}\label{d-cyclic}
 Let $n$, $n'$ be positive integers and $F$ be a finite field with
characteristic coprime to both $n$ and $n'$.
Let $C_{g,g'}$ be the double cyclic code generated by
$\big(g(X),g'(X)\big)$ where $g(X)=\sum_{j=0}^{n-1}g_jX^j$
and $g'(X)=\sum_{j'=0}^{n'-1}g'_{j'}X^{j'}$.
Let $m={\rm lcm}(n,n')$ which denotes
the least common multiple of $n$ and $n'$. Let $\ell=\gcd(n,n')$,
$$\textstyle
 d=\deg\frac{\gcd\!\big(g(X),\,X^n-1\big)\cdot
     \gcd\!\big(g'(X),\,X^{n'}-1\big)\cdot(X^\ell -1)}
   {\gcd\!\big(g(X)g'(X),\,X^\ell-1\big)}
$$
and $r=n+n'-d$. Then the dimension
$\dim C_{g,g'}=r$ and the first $r$ rows of the double circulant
matrix $M(g,g')_{m\times(n+n')}$ associated with $g(X)$ and $g'(X)$
form a generator matrix $G$ of $C_{g,g'}$.
\end{theorem}

\begin{example}\rm
Let $F={\Bbb Z}_5$ (the integer residue ring modulo $5$)
be the finite field with $5$ elements,
$n=2$, $n'=3$, $g(X)=-1+X$, $g'(X)=-2+X+X^2$.
Then $\ell=\gcd(n,n')=1$, $\gcd\big(g(X), X^2-1\big)=X-1$,
$\gcd\big(g'(X), X^3-1\big)=X-1$,
$\gcd\big(g(X)g'(X), X-1\big)=X-1$.
By Theorem \ref{d-cyclic}, $\dim C_{g,g'}=3$.
Since the double circulant matrix
$$
M(g,g')_{(2+3)\times(2+3)}=
\begin{pmatrix}
  -1&1&-2&1&1\\ 1&-1&1&-2&1\\-1&1&1&1&-2\\
  1&-1&-2&1&1\\ -1&1&1&-2&1
\end{pmatrix}.
$$
we get a generator matrix $G$ of $C_{g,g'}$ as follows:
$$
G=
\begin{pmatrix}
  -1&1&-2&1&1\\ 1&-1&1&-2&1\\-1&1&1&1&-2
\end{pmatrix}.
$$
\end{example}

\section{Extensions and questions}\label{multiple circulant}

A natural generalization of double circulant matrices
in Eqn \eqref{GDC} is as follows.

\begin{definition}\label{m-circulant}\rm
Let $F$ be a field, $m$, $n_1,\cdots,n_k$ be positive integers
such that every $n_i$ is coprime to ${\rm char}\,F$,
and $g_i(X)=\sum_{j=0}^{n_i-1}g_{ij}X^j\in F[X]$ for $i=1,\cdots,k$.
Let $M(g_i)_{m\times n_i}$ for $i=1,\cdots,k$
be the generalized circulant matrices defined in Eqn \eqref{GCM}.
The concatenation side-by-side of the $k$ matrices $M(g_i)_{m\times n_i}$
for $i=1,\cdots,k$:
$$
M(g_1,\cdots,g_k)_{m\times(n_1+\cdots+n_k)}
=\Big(M(g_1)_{m\times n_1}\,\big|\,\cdots\,\big|\,M(g_k)_{m\times n_k}\Big)
$$
is called a {\em multiple circulant matrix} associated with the polynomials
$g_i(X)$ for $i=1,\cdots,k$.
\end{definition}

Theorem \ref{thm d-circulant} and
Corollary \ref{cor thm for d-circulant}
are easily extended. To state them,
we show another description of the quantity $n+n'-d$ in
Theorem \ref{thm d-circulant}.

\begin{remark}\rm
Recall that $\ell=\gcd(n,n')$ in Section 3, and
by checking roots of polynomials, we showed that
$\gcd\big(X^n-1,\,X^{n'}-1\big)=X^{\ell}-1$, see Eqn \eqref{gcd n n'}.
By a similar argument, we can show that, for any $g(X),g'(X)\in F[X]$,
\begin{equation}\label{dd=}\textstyle
\begin{array}{l}
{\rm lcm}\Big(\frac{X^n-1}{\gcd\!\big(g(X),X^n-1\big)},~
   \frac{X^{n'}-1}{\gcd\!\big(g'(X),X^{n'}-1\big)}\Big)\\
~=\frac{(X^n-1)(X^{n'}-1)}
{\gcd\!\big(g(X),X^{n}-1\big)\cdot\gcd\!\big(g'(X),X^{n'}-1\big)}
       \cdot\frac{\gcd\!\big( g(X)g'(X),~X^\ell-1\big)}{X^\ell -1},
  \end{array}
\end{equation}
where ${\rm lcm}$(-,-) denotes the least common multiple. Thus
$$\textstyle
n+n'-d=\deg{\rm lcm}\Big(\frac{X^n-1}{\gcd(g(X),X^n-1)},~
   \frac{X^{n'}-1}{\gcd(g'(X),X^{n'}-1)}\Big),
$$
where $d$ is defined in Eqn \eqref{d=}.
\end{remark}

\begin{theorem}\label{thm m-circulant}
Let notations be as in Definition \ref{m-circulant}. Let
$$\textstyle
 s=\deg{\rm lcm}\Big(\frac{X^{n_1}-1}{\gcd\!\big(g_1(X),X^{n_1}-1\big)},
 ~\cdots,~\frac{X^{n_k}-1}{\gcd\!\big(g_k(X),X^{n_k}-1\big)}\Big),
$$
and let $r=\min\{m,s\}$. Then the rank
$$
{\rm rank}\,M(g_1,\cdots,g_k)_{m\times(n_1+\cdots+n_k)}=r,
$$
and the first consecutive $r$ rows of
$M(g_1,\cdots,g_k)_{m\times(n_1+\cdots+n_k)}$
are linearly independent.
\end{theorem}

\pf
The proof is similar to what we did in Section \ref{d-c matrix}, we sketch it.
Take primitive $n_i$-th root $\omega_i$ of unity for $i=1,\cdots,k$.
Similarly to Eqn \eqref{FTM-d}, define a matrix:
$$
\Phi(\omega_1,\cdots,\omega_k)=\Big(\Phi(\omega_1)_{m\times n_1}
 \,\big|\,\cdots\,\big|\,\Phi(\omega_k)_{m\times n_k}\Big).
$$
Similarly to Eqn \eqref{MD=}, denote
\begin{equation*}
\begin{array}{l}
 D(\omega_1,\cdots,\omega_k)={\rm diag}\Big(\Phi(\omega_1)_{n_1\times n_1},
  \,\cdots,\,\Phi(\omega_k)_{n_k\times n_k}\Big),\\
 D(g_1,\cdots,g_k)={\rm diag}\big(g_1(\omega_1^0),\cdots,g_1(\omega_1^{n_1 -1}),
  \,\cdots,\,
   g_k(\omega_k^0),\cdots,g_k(\omega_k^{n_k -1})\big).
\end{array}
\end{equation*}
Then, similarly to Lemma \ref{C-F-d}, we have
\begin{equation*}
M(g_1,\cdots,g_k)_{m\times(n_1+\cdots+n_k)}\!\cdot\! D(\omega_1,\cdots,\omega_k)
 =\Phi(\omega_1,\cdots,\omega_k)\!\cdot\! D(g_1,\cdots,g_k).
\end{equation*}

As we said in the proof of Theorem \ref{thm d-circulant},
the next key step is to determine which of the columns of the matrix
$\Phi(\omega_1,\cdots,\omega_k)\!\cdot\! D(g_1,\cdots,g_k)$
contribute essentially to the rank of the matrix.
They are determined by two observations.

$\bullet$~
Each diagonal element of $D(g_1,\cdots,g_k)$ is of the form
$g_i(\omega_i^j)$. If $g_i(\omega_i^j)=0$ then the corresponding column of
$\Phi(\omega_1,\cdots,\omega_k)\!\cdot\! D(g_1,\cdots,g_k)$ is zero.
Hence only such columns remain to be considered,
which correspond to $g_i(\omega_i^j)\ne 0$, i.e.,
$\omega_i^j$ is a root of $\frac{X^{n_i}-1}{\gcd\!\big(g_i(X),X^{n_i}-1\big)}$.

$\bullet$~
Each column of $\Phi(\omega_1,\cdots,\omega_k)$ is of the form
$V(\omega_i^j)_m$. If $\omega_{i}^{j}=\omega_{i'}^{j'}$,
then $V(\omega_i^j)_m=V(\omega_{i'}^{j'})_m$ which contribute
at most one to the rank of the matrix
$\Phi(\omega_1,\cdots,\omega_k)\!\cdot\! D(g_1,\cdots,g_k)$.

Therefore, corresponding to each root of the least common multiple
$$\textstyle
 {\rm lcm}\left(\frac{X^{n_1}-1}{\gcd\!\big(g_1(X),X^{n_1}-1\big)},
 ~\cdots,~\frac{X^{n_k}-1}{\gcd\!\big(g_k(X),X^{n_k}-1\big)}\right),
$$
we take a column of $\Phi(\omega_1,\cdots,\omega_k)$ and a column
of $D(g_1,\cdots,g_k)$. Then, similarly to Eqn\eqref{Phi D},
we get an $m\times s$ matrices $\Phi$ and
an $s\times s$ diagonal matrix $D$ such that
$$
{\rm rank}(\Phi\cdot D)=
{\rm rank}\big(\Phi(\omega_1,\cdots,\omega_k)\!\cdot\! D(g_1,\cdots,g_k)\big).
$$
Finally, the proof can be finished in the same way as
the last paragraph of the proof of
Theorem \ref{thm d-circulant}.
\qed

In the same way as for Corollary \ref{cor thm for d-circulant},
we have a corollay.

\begin{corollary}\label{cor m-circulant}
Any consecutive $r$ rows of $M(g_1,\cdots,g_k)_{m\times(n_1+\cdots+n_k)}$
are linearly independent.
\end{corollary}

Theorem \ref{thm m-circulant} and Corollary \ref{cor m-circulant}
can be used to construct generator matrices of a kind of
multiple cyclic codes (or generalized cyclic codes), see \cite{EY}.

\medskip
We conclude the paper with few questions on double circulant matrices.

Though the double circulant matrix $M(g,g')_{m\times(n+n')}$
defined in Eqn \eqref{GDC} is not square in general,
we can consider the square case $M(g,g')_{(n+n')\times(n+n')}$
when we take $m=n+n'$.
Then a fundamental question comes up naturally.

\medskip\noindent
{\bf Question 1.}~ {\em How to find eigenvalues and eigenvectors of
the square double circulant matrix
$M(g,g')_{(n+n')\times(n+n')}$? }

\medskip
From Corollary \ref{eigen 0}, we see that
$0$ is an eigenvalue of $M(g,g')_{(n+n')\times(n+n')}$,
and $d$ eigenvectors of the eigenvalue $0$
which are linearly independent are obtained.
However, we didn't show that the multiplicity of the eigenvalue $0$
is equal to $d$.

\medskip\noindent
{\bf Question 2.}~ {\em Is the square double circulant matrix
$M(g,g')_{(n+n')\times(n+n')}$ diagonalizable? }

\section*{Acknowledgements}
The research of the authors is supported by NSFC
with grant numbers 11271005.

\end{document}